\documentclass{amsart}

\usepackage{amssymb}
\usepackage{amsfonts}
\usepackage{amsmath}
\usepackage{amsthm}
\usepackage[pdftex,bookmarks,colorlinks]{hyperref}
\usepackage[pdftex]{hyperref}
\usepackage{url}
\usepackage{color}
\usepackage{graphicx}
\usepackage{soul}
\usepackage[small,nohug,heads=LaTeX]{diagrams}

\newtheorem{theorem}{Theorem}[section]

\newtheorem{corollary}[theorem]{Corollary}

\newcommand{\bq}{\begin{quote}}
\newcommand{\eq}{\end{quote}}
\newcommand{\bi}{\begin{itemize}}
\newcommand{\ei}{\end{itemize}}
\newcommand{\bd}{\begin{description}}
\newcommand{\ed}{\end{description}}
\newcommand{\ben}{\begin{enumerate}}
\newcommand{\een}{\end{enumerate}}
\newcommand{\bbm}{\begin{bmatrix}}
\newcommand{\ebm}{\end{bmatrix}}
\newcommand{\bea}{\begin{eqnarray*}}
\newcommand{\eea}{\end{eqnarray*}}

\def\bF{\mathbb{F}}

\def\GF{\mathbb{F}}

\def\KK{\mathsf{K}}

\def\PP{\mathsf{P}}
\def\RR{\mathsf{R}}

\def\XX{\mathsf{X}}
\def\YY{\mathsf{Y}}
\def\ZZ{\mathsf{Z}}

\newcommand{\bx}[1]{\mathsf{#1}}     
\newcommand{\bxg}[1]{\boldsymbol #1} 

\def\Aut{\mathop{{\rm Aut}}}

\def\2G2{\ensuremath{^2{\rm G}_2}}
\def\sl{{\rm{SL}}}

\def\psl{{\rm{PSL}}}

\def\so{{\rm{SO}}}

\def\ppsl{ ( {\rm{P}} ) {\rm{SL}}}

\definecolor{darkgreen}{rgb}{0,0.6,0}

\newcommand{\encr}{\ensuremath{\vDash}}

\AtBeginDocument{%
  \def\MR#1{}
  }

\begin{document}

\title{Homomorphic encryption and some black box attacks\footnote{It will appear in Lecture Notes in Computer Science (LNCS).}}

\author{Alexandre Borovik}
\address{School of Mathematics, University of Manchester, UK; alexandre $\gg {\rm at} \ll $ borovik.net}
\author{\c{S}\"{u}kr\"{u} Yal\c{c}\i nkaya}
\address{Department of Mathematics, Istanbul University, Turkey; sukru.yalcinkaya $\gg {\rm at} \ll $ istanbul.edu.tr}

\maketitle

\begin{abstract}
This paper is a compressed summary of some principal definitions and concepts in the approach to the black box algebra being developed by the authors \cite{BY2020-sl2,BY2014,BY2018}.    We suggest that black box algebra could be useful in cryptanalysis of homomorphic encryption schemes \cite{fontaine09.41}, and that homomorphic encryption is an area of research where cryptography and black box algebra may benefit from exchange of ideas.
\end{abstract}

\section{Homomorphic encryption}
\label{sec:homomorphic}

``Cloud computing'' appears to be a hot topic in information technology; in a nutshell, this is the ability of small and computationally weak  devices to delegate hard resource-intensive computations to third party (and therefore untrusted) computers. To ensure the privacy of the data, the untrusted computer should receive data in an encrypted form but still being able to process it. It means that encryption should preserve algebraic structural properties of the data. This is one of the reasons for popularity of the idea of \emph{homomorphic encryption} \cite{acar2018,aguilar2013recent,DBLP:journals/corr/DyerDX17,fontaine09.41,Gentry:2009:FHE:1536414.1536440,gentry11.129,DBLP:journals/corr/PrasannaA15,DBLP:journals/corr/Rass13,DBLP:journals/corr/abs-1305-5886,DBLP:journals/corr/Sharma13a,DBLP:journals/corr/TebaaH14}
which we describe here with some simplifications aimed at clarifying connections with black box algebra (as defined in Section \ref{sec:BB1--BB3}).

\subsection{Homomorphic encryption: basic definitions}\label{subsec:hom_enc}

Let $A$ and $\XX$ denote the sets of plaintexts and ciphertexts, respectively, and assume that we have some (say, binary) operators  $\boxdot_A$ on $A$ needed for processing data and corresponding  operators $\boxdot_\XX$  on $\XX$. An encryption  function $E$ is  \emph{homomorphic} if
\[
E(a_1 \boxdot_A a_2) = E(a_1) \boxdot_\XX E(a_2)
\]
for all plaintexts $a_1$, $a_2$ and all operators on $A$.

Suppose that  Alice is the owner of data represented by plaintexts in $A$ which she would like to process using operators $\boxdot_A$ but has insufficient computational resources, while Bob has computational facilities for processing ciphertexts using operators $\boxdot_\XX$. Alice may wish to enter into a contract with Bob; in a realistic scenario,  Alice is one of the many customers of the encrypted data processing service run by Bob, and all customers use the same ambient structure $A$ upto isomorphism and formats of data and operators which are for that reason are likely to be known to Bob. What is not known to Bob is the specific password protected encryption used by Alice. This is what is known in cryptology as Kerckhoff's Principle: \emph{obscurity is no security}, the security of encryption should not rely on details of the protocol being held secret; see \cite{fontaine09.41} for historic details.

Alice encrypts plaintexts $a_1$ and $a_2$ and sends ciphertexts $E(a_1)$ and $E(a_2)$ to Bob, who computes
\[
\bx{x} = E(a_1) \boxdot_\XX E(a_2)
\]
without having access to the content of plaintexts $a_1$ and $a_2$, then return the output $\bx{x}$ to Alice who decrypts it using the decryption function $E^{-1}$:
\[
E^{-1}(\bx{x}) = a_1 \boxdot_A a_2
\]
In this set-up, we say that the homomorphic encryption scheme is \emph{based} on the algebraic structure $A$ or the homomorphism $E$ is a \emph{homomorphic encryption of} the algebraic structure $A$.

To simplify exposition, we assume that the encryption function $E$ is deterministic, that is, $E$  establishes a one-to-one correspondence between $A$ and $\XX$. Of course, this is a strong assumption in the cryptographic context; it is largely unnecessary for our analysis, but, for the purposes of this paper, allows us to avoid technical details and makes it easier to explain links with the black box algebra.

\subsection{Back to algebra}

In algebraic terms, $A$ and  $\XX$ as introduced above are algebraic structures with operations on them which we refer to as \emph{algebraic operations} and
\(
E: A \longrightarrow \XX
\)
is a homomorphism. In this paper we assume that the algebraic structure $A$ is finite as a set. This is not really essential for our analysis, many observations are relevant for the infinite case as well, but handling probability distributions (that is, random elements) on infinite sets is beyond the scope of the present paper.

We discuss a class of potential attacks on homomorphic encryption of $A$. Our discussion is based on a simple but fundamental fact of algebra that a map
\(
E:A \longrightarrow \XX
\)
of algebraic structures of the same type is a homomorphism if and only if its graph
\[
\Gamma(E) = \{(a,E(a)) \mid a\in A\}
\]
is a substructure of $A\times \XX$, that is, closed under all algebraic operations on $A\times \XX$. Obviously, $\Gamma(E)$ is isomorphic to $A$ and we shall note the following observation:
\begin{quote}
\textbf{if an algebraic structure $A$ has a rich internal configuration (has many substructures with complex interactions between them), the graph $\Gamma(E)$ of a homomorphic encryption $E: A \longrightarrow \XX$ also has a rich (admittedly hidden) internal configuration, and this could make it vulnerable to an attack from Bob.}
\end{quote}

We suggest that
\bq
\textbf{before attempting to develop a homomorphic encryption scheme based on a particular algebraic structure $A$, the latter needs to be examined by black box theory methods -- as examples in this paper show, it could happen that \emph{all} homomorphic encryption schemes on $A$ are insecure.}
\eq

\section{Black box algebra}

\subsection{Axiomatic description of black box algebraic structures}
\label{sec:BB1--BB3}

A \emph{black box algebraic structure} $\XX$ is a black box (device, algorithm, or oracle) which produces and  operates with $0$--$1$ strings of uniform length $l(\XX)$ encrypting (not necessarily in a unique way) elements of some fixed algebraic structure $A$: if $\bx{x}$ is one of these strings  then it corresponds to a unique (but unknown to us) element $\pi(\bx{x}) \in A$.
Here, $\pi$ is the decrypting map, not necessarily known to us in advance.   We call the strings produced or computed by $\XX$ \emph{cryptoelements}.

Our axioms for black boxes are the same as in \cite{BY2020-sl2,BY2014,BY2018}, but stated in a more formal language.
\bi
\item[\textbf{BB1}] On request, $\XX$ produces a `random' cryptoelement $\bx{x}$ as a string of fixed length $l(\XX)$, which depends on $\XX$, which encrypts an element $\pi(x)$ of some fixed explicitly given  algebraic structure $A$; this is done in  time polynomial in $l(\XX)$. When this procedure is repeated, the  elements $\pi(\bx{x}_1), \pi(\bx{x}_2),\dots$ are independent and uniformly distributed in $A$.
\ei

To avoid messy notation, we assume that operations on $A$ are unary or binary; a general case can be treated in exactly the same way.

\bi
\item[\textbf{BB2}] On request,  $\XX$ performs algebraic operations on the encrypted strings which correspond to operations in  $A$ in a way which makes the map $\pi$ (unknown to us!) a homomorphism: for every binary (unary case is similar) operation $\boxdot $ and strings $\bx{x}$ and $\bx{y}$ produced or computed by $\XX$,
\[
\pi(\bx{x} \boxdot \bx{y}) = \pi(\bx{x}) \boxdot \pi(\bx{y}).
\]
\ei

It should be noted that we do not assume the existence of an algorithm which allows us to decide whether a specific string can be potentially produced by $\XX$; requests for operations on strings can be made only in relation to cryptoelements previously output by $\XX$. Also, we do not make any assumptions on probabilistic distribution of cryptoelements.

\begin{itemize}
\item[\textbf{BB3}] On request, $\XX$ determines, in time polynomial in $l(\XX)$, whether two cryptoelements $\bx{x}$ and $\bx{y}$ encrypt the same element in $A$, that is, check whether $\pi(\bx{x}) = \pi(\bx{y})$.
\end{itemize}

We say in this situation that  a black box $\XX$ \emph{encrypts} the algebraic structure $A$ and we denote this as $\XX \encr A$.

Clearly, in black box problems, the decrypting map $\pi$ is not given in advance. However, it is useful to think about any algebraic structure (say, a finite field) implemented on a computer as a trivial black box, with $\pi$ being the identity map, and with random elements produced with the help of a random number generator. In this situation, obviously, the axioms BB1 -- BB3 hold.

In our algorithms, we have to build new black boxes from existing ones and work with several black box structures at once: this is why we have to keep track of the length $l(\bx{\XX})$ on which a specific black box $\XX$ operates. For example, it turns out in \cite{BY2018} that it is useful to consider an automorphism of $A$ as a graph in $A\times A$. This produces an another algebraic structure isomorphic to $A$ which can be seen as being encrypted by a black box $\ZZ$ producing, and operating on, certain pairs of strings from $\XX$, see \cite{BY2018} for more examples. In this case, clearly, $l(\ZZ)  =2l(\XX)$.

\subsection{Morphisms}

Given two  black boxes $\XX$ and $\YY$ encrypting  algebraic structures $A$ and  $B$, respectively, we say that a map $\bxg{\phi}$ which assigns
strings produced by $\XX$  to strings produced by $\YY$  is a \emph{morphism} of black boxes, if
\bi
\item the map $\bxg{\phi}$ is computable in time polynomial in $l(\XX)$ and $l(\YY)$, and
\item there is a homomorphism $\phi:A \to B$ such that the following diagram  is commutative:
\begin{diagram}
\XX &\rTo^{\bxg{\phi}} &\YY\\
\dDotsto_{\pi_{\XX}} & &\dDotsto_{\pi_{\YY}}\\
A &\rTo^{\phi} & B
\end{diagram}
where $\pi_\XX$ and $\pi_\YY$ are the canonical projections of $\XX$ and $\YY$ onto $A$ and $B$, respectively.
\ei
We  say in this situation that a morphism $\bxg{\phi}$ \emph{encrypts} the homomorphism $\phi$ and call $\bxg{\phi}$ bijective, injective, etc., if $\phi$ has these properties.

\subsection{Construction and interpretation}

Construction of a new black box $\YY$ in a given black box $\XX\encr A$ can be formally described as follows.

Strings of $\YY$ are concatenated  $n$-tuples of strings $(\bx{x}_1,\dots, \bx{x}_n)$ from $\XX$ produced by a polynomial time algorithm which uses operations on $\XX$; new operations on $\YY$ are also polynomial time algorithms running on $\XX$, as well as the algorithm for checking the new identity relation $=_{\YY}$ on $\YY$.

If this is done in a consistent way and axioms BB1--BB3 hold in $\YY$, then $\YY$ encrypts an algebraic structure $B$ which can be obtained from the structure $A$ by a similar construction, with algorithms replaced by description of their outputs by formulae of first order language in the signature of $A$. At this point we are entering the domain of model theory, and full discussion of this connection can be found in our forthcoming paper \cite{BY2020-models}. Here we notice only that in model theory $B$ is said to be \emph{interpreted} in $A$, and if $A$ is in its turn  interpreted in $B$ then $A$ and $B$ are called \emph{bi-intrepretable}. A recent result on bi-interpretability between Chevalley groups and rings, relevant to our project is \cite{Segal-Tent2020}.

\subsection{A few historic remarks}

Black box algebraic structures had been introduced by Babai and Szem\'{e}redi \cite{babai84.229} in the special case of groups as an idealized setting for randomized algorithms for solving permutation and matrix group problems in computational group theory.  Our Axioms BB1--BB3 are a slight modification -- and generalization to arbitrary algebraic structures -- of their original axioms.

So far, it appears that only finite groups, fields, rings, and, very recently, projective planes (in our paper \cite{BY2018}) got a black box treatment. In the case of finite fields, the concept of a black box field can be traced back to Lenstra Jr \cite{lenstra91.329} and Boneh and Lipton \cite{Boneh96.283}, and in the case of rings -- to Arvind \cite{arvind06.126}.

A higher level of abstraction introduced in our papers produces new tools allowing us to solve problems which previously were deemed to be intractable. For example, recently, a fundamental problem of constructing a unipotent element in black box groups encrypting $\psl_2$ was solved in odd characteristics via constructing a black box projective plane and its underlying black box field \cite{BY2018}. There is an analogous recognition algorithm for the black box groups encrypting $\psl_2$ in even characteristic \cite{kantor15.16}.

\subsection{Recognition of black box fields}
\label{sec:fields}

A \emph{black box} (finite) \emph{field} $\KK$ is a black box operating on $0$-$1$ strings of uniform length which encrypts some finite field $\GF$. The oracle can compute $\bx{x}+\bx{y}$, $\bx{x}\bx{y}$, and $\bx{x}^{-1}$ (the latter for $\bx{x} \ne 0$) and decide whether $\bx{x}=\bx{y}$ for any strings $\bx{x},\bx{y} \in \KK$. Notice in this definition that the characteristic of the field is not known. Such a definition is needed in our paper \cite{BY2018} to produce black box group algorithms which does not use characteristic of the underlying field. If the characteristic $p$ of $\KK$ is known then we say that $\KK$ is a \emph{black box field of known characteristic $p$}. We refer the reader to \cite{Boneh96.283,maurer07.427} for more details on black box fields of known characteristic and their applications to cryptography.

The following theorem is a reformulation of the fundamental results in \cite{maurer07.427}.

\begin{theorem} \label{th:bbfields}
Let\/ $\KK \encr \mathbb{F}_{p^n}$  be a black box field of known characteristic $p$ and $\KK_0$ the  prime subfield of\/ $\KK$. Then the problem of finding two way morphisms between $\KK$ and\/ $\mathbb{F}_{p^n}$ can be reduced to the same problem for $\KK_0$ and\/ $\mathbb{F}_p$. In particular,
\bi
\item  a  morphism $\KK_0 \longrightarrow \mathbb{F}_p$ can be extended in time polynomial in the input length\/ $l(\KK)$ to a  morphism $\KK \longrightarrow \mathbb{F}_{p^n};$
\item there is a  morphism  $\mathbb{F}_{p^n} \longrightarrow \KK$ computable in time polynomial in  $l(\KK)$.
    \ei
\end{theorem}

Here and in the rest of the paper, ``efficient'' means ``computable in time polynomial in the input length''.

In our terminology (Section~\ref{sec:structural_proxy}), Theorem \ref{th:bbfields}  provides a \emph{structural proxy} for black box fields of known characteristic. Indeed, if $\KK$ is a black box field of known characteristic $p$, then we can construct an isomorphism $\GF_p =\mathbb{Z}/p\mathbb{Z} \longrightarrow \KK_0$ by the map
\[
m \mapsto \bx{1}+\bx{1}+\cdots+ \bx{1} \quad (m \mbox{ times})
\]
where $\bx{1}$ is the unit in $\KK_0$; it is computable in linear in $\log p$ time by double-and-add method. We say that $p$ is \emph{small} if it is computationally feasible to make a lookup table for the inverse $\KK_0 \longrightarrow \GF_p$ of this map. Construction of a morphism $\GF_p \longleftarrow \KK_0$ remains an open problem. However, we can observe that

\begin{corollary}\label{cor:bbfields} Let $\KK \encr \bF_{p^n}$, where $p$ is a known small prime number. Then there exist two way morphisms between $\KK$ and $\bF_{p^n}$.
\end{corollary}

\subsection{Construction of a structural proxy}\label{sec:structural_proxy}

Most groups of Lie type (we exclude $^2B_2$, $^2F_4$ and $^2G_2$ to avoid technical details) can be seen as functors $G: \mathcal{F} \longrightarrow \mathcal{G}$ from the category of fields $\mathcal{F}$ with an automorphism of order $\leqslant 2$ to the category of groups $\mathcal{G}$. There are also other algebraic structures which can be defined in a similar way as functors from $\mathcal{F}$, for example projective planes or simple Lie algebras (viewed as rings). The following problem is natural and, as our results show, useful in this context.

\bi
\item[] \emph{Construction of a structural proxy:}  Suppose that we are given  a black box structure $\XX\encr A(\bF)$. Construct, in time polynomial  in $l(\XX)$,
    \bi
     \item a black box field $\KK\encr \bF$, and
     \item  two way bijective morphisms $A(\KK) \longleftrightarrow \XX$.
     \ei
\ei

If we construct a black box field $\KK$ by using $\XX$ as a computational engine, then we can construct the natural representation $A(\KK)$ of the structure $A$ over the black box field $\KK$. By Theorem \ref{th:bbfields}, we can construct a polynomial time isomorphism $\bF_q \longrightarrow \KK$ which further provides an isomorphism $A(\bF_q) \longrightarrow A(\KK)$  completing \emph{a structure recovery} of $\XX$.

Structural proxies and structure recovery play a crucial role for algorithms developed in Theorem \ref{homomorphic-SL2}.
We summarize relevant results about constructing structural proxies of black box algebraic structures from our papers \cite{BY2020-sl2,BY2018}.

\begin{theorem}\label{theorem:proxies}
We can construct structural proxies for the following black box structures.

\bi
\item[{\rm (a)}] $\PP \encr \mathbb{P}^2(\bF)$,  a projective plane with a polarity encrypting a projective plane  $\mathbb{P}^2(\bF)$ over a finite field $\bF$ of odd characteristic.

\item[{\rm (b)}] $\XX \encr \so_3(\bF), \ppsl_2(\bF)$ over a finite field $\bF$ of unknown odd characteristic, under the assumption that we know a global exponent $E$ of $\XX$, that is, $E$ such that $\bx{x}^E = \bx{1}$ for all $\bx{x}\in \XX$ and $\log E$ is polynomially bounded in terms of $l(\XX)$.

\item[{\rm (c)}] $\RR \encr {\rm M}_{2\times 2}(\bF_q)$, a black box ring encrypting the ring of $2\times 2$ matrices over the known finite field $\bF_q$ of odd characteristic.
\ei
\end{theorem}

\subsection{Black boxes associated with homomorphic encryption}

As explained in Subsection~\ref{subsec:hom_enc}, we assume that the algebraic structure $A$ of plaintexts is represented in some standard form known to Bob. In agreement with the standard language of algebra -- and with our terminology in \cite{BY2018} -- we shall use the words \emph{plain element} or just \emph{element} in place of  `plaintext' and \emph{cryptoelement}  in place of `ciphertext'.

Let $A$ be a set of plain elements, $\XX$ a set of cryptoelements, and  $E$ be the encryption function, that is, an isomorphism $E: A \longrightarrow  \XX$.

Supply of random cryptoelements from $\XX$ postulated in Axiom BB1 can be achieved by sampling a big dataset of cryptoelements provided by Alice, or computed on request from Alice. The computer system controlled by Bob performs algebraic operations referred to in Axiom BB2.

Axiom BB3 is redundant under the assumption that $E: A \longrightarrow \XX$ is a bijection but it gives us more freedom to construct new black boxes, for example, homomorphic images of $\XX$. Axiom BB3 could also be useful for handling another quite possible scenario: For Alice, the cost of computing homomorphisms $E$ and $E^{-1}$ could be higher than the price charged by Bob for processing cryptoelements. In that case, it could be cheaper to transfer initial data to Bob (in encrypted form) and ask Bob to run a computer programme which uses the black box but does not send intermediate values back to Alice, returns only the final result; checking equality of cryptoelements becomes unavoidable.

\section{A black box attack on homomorphic encryption}
\label{sec:impersonalisation}

We assume that Bob can accumulate a big dataset of cryptoelements sent from/to Alice, or intermediate results from running Alice's programme, and that he can feed, without Alice's knowledge, cryptoelements into a computer system (the \emph{black box}) which performs operations on them, and retain the outputs for peruse -- again without Alice's knowledge.  Bob's aim is to compute the decryption function $E^{-1}$ efficiently, that is, in time polynomial in terms of the lengths of plain elements and cryptoelements involved.

\subsection{Bob's attack}
\label{Bobs-attack}

As we discussed in Section \ref{subsec:hom_enc}, we can assume that Bob knows the algebraic structure $A$. Bob's aim is to find an efficient algorithm which maps cryptoelements from $\XX$ to elements in $A$ and vice versa while preserving the algebraic operations on $\XX$ and $A$. This means solving the \emph{constructive recognition problem} for $\XX$,
that is, finding bijective morphisms
\[
\alpha: \XX \longrightarrow A \mbox{ and }
\beta: A  \longrightarrow \XX
\]
such that $\alpha\circ\beta$ is the identity map on $A$.

Assume that Bob solved the constructive recognition problem and can efficiently compute $\alpha$ and $\beta$.

Alice's encryption function is a map $E: A \longrightarrow \XX$; the composition $\delta = \alpha \circ E$ is an automorphism of $A$. Therefore Bob reads not Alice's plaintexts $a\in A$, but their images $\delta(a) = \alpha(E(a))$ under an automorphism $\delta$ of $A$ still unknown to him. This means that
\bq
\textbf{solving the constructive recognition problem for $\XX$ reduces the problem of inverting the encryption homomorphism $E: A \longrightarrow \XX$ to a much simpler problem of inverting the automorphism}
\[
\delta: A \longrightarrow A.
\]
\eq
We are again in the situation of homomorphic encryption, but this time the sets of plaintexts and ciphertexts are the same. One would expect that this encryption is easier to break. For example, if Bob can guess the plaintexts of a few cryptoelements, and if the automorphism group $\Aut A$ of $A$ is well understood, computation of $\delta$ and $\delta^{-1}$ could be a more accessible problem than the constructive recognition for $\XX$. For example, automorphism groups of finite fields are very small, and in that case $\delta^{-1}$ can be found by direct inspection.

{As soon as $\delta^{-1}$ is known, Bob knows
$
E^{-1} = \delta^{-1}\circ\alpha
$
and can decrypt everything. Moreover, since
$
E =\beta \circ \delta,
$
the map $E$ is also known and allows Bob to return to Alice cryptoelements which encrypt plaintexts of Bob's choice.

We suggest that this approach to analysis of homomorphic encryption is useful because it opens up connections to black box algebra. Indeed the theory of black box structures is reasonably well developed  for groups and fields, and its methods could provide insight into assessment of security of other algebraic structures if any are proposed for use in homomorphic encryption.

\section{Application of Theorem \ref{theorem:proxies} to homomorphic encryption}
\label{sec:two-groups}

The procedures described in Theorem \ref{homomorphic-SL2} below are reformulations of the principal results of our
Theorem \ref{theorem:proxies} in a homomorphic encryption setup. They demonstrate the depth of structural analysis involved
and suggest that a similarly deep but revealing structural theory can be developed for other algebraic structures if they are sufficiently rich (`rich' here can mean, for example, `bi-interpretable with a finite field'). Also, it is worth noting that the procedures do not use any assumptions about the encryption homomorphism $E$, the analysis is purely algebraic.

\begin{theorem}
\label{homomorphic-SL2}
Assume that Alice and Bob run a homomorphic encryption protocol over the group $A = \sl_2(\bF_q)$, $q$ odd, with Bob doing computations with cryptoelements using a black box $\XX\encr A$. Assume that Bob  knows $A$, including the representation of the field $\bF_q$ used by Alice. Then, by Theorem \emph{\ref{theorem:proxies}}, Bob can construct a structural proxy $\sl_2(\KK)$ for $\XX$. Moreover:
\bi
\item[{\rm (a)}] If, in addition, Bob has  two way bijective morphisms between a black box field $\KK$ and an explicitly given field $\bF_q$ \emph{(}see Corollary \emph{\ref{cor:bbfields})}, he gets  two way bijective morphisms $\XX \longleftrightarrow \sl_2(\bF_q)$.

\item[{\rm (b)}] Under assumptions of \emph{(a)}, Bob gets an image of Alice's data transformed by an automorphism $\delta: \sl_2(\bF_q)\longrightarrow \sl_2(\bF_q)$ since Alice's group $A$ is an explicitly given $\sl_2(\bF_q)$.

\item[{\rm (c)}]  Automorphisms  of the group $\sl_2(\bF_q)$ are well known: every automorphism is a product of an inner automorphism  and a field automorphism induced by an automorphism of the field $\bF_q$.  Therefore if  Bob can run a few instances of known plaintexts attacks against Alice, he can  compute the automorphism $\delta$ and after that read plaintexts of all  Alice's cryptoelements.

\item[{\rm (d)}] Moreover, under assumptions of \emph{(a)} and \emph{(c)},  Bob can compute the inverse of $\delta$ and pass to Alice, as answers to Alice's requests, values of his choice.
\ei

\end{theorem}

Items (c) and (d) in Theorem \ref{homomorphic-SL2} look as serious vulnerabilities of homomorphic encryptions of the groups $\sl_2(\bF_q)$. We conclude that homomorphic encryption of  groups $\sl_2(\bF_q)$ is no more secure  than homomorphic encryption of the field $\bF_q$. As a consequence of Theorem \ref{th:bbfields}, homomorphic encryption of $\sl_2(\bF_q), q = p^k,$ does not survive a known plaintext attack when the prime  $p >2$ is small.

We think that this is a manifestation of a more general issue:  for small odd primes $p$, there are no secure homomorphic encryption schemes based on sufficiently rich (say, bi-interpretable with finite fields) algebraic structures functorially defined over finite fields of characteristic $p$.

\section*{Acknowledgements}

The authors worked on this paper during their visits to the Nesin Mathematics Village, Turkey. We thank Jeff Burdges, Adrien Deloro, Alexander Konovalov, and Chris Stephenson for fruitful advice, and the referees for their most perceptive comments.

\end{document}